\documentclass[twoside]{ima}
  \newif\ifpdf\ifx\pdfoutput\undefined\pdffalse\else\pdfoutput=1\pdftrue\fi

\usepackage{amssymb}

\usepackage{amsmath}
\usepackage{float}
\usepackage{graphicx}

\floatstyle{plain}
\restylefloat{figure}

\DeclareMathOperator{\ns}{ns}
\DeclareMathOperator{\clos}{Cl}

\begin{document}

\title{POLARS OF REAL SINGULAR PLANE CURVES}
\author{HEIDI CAMILLA MORK
\and RAGNI PIENE\thanks{CMA, Department of Mathematics, University of Oslo, P.
O. Box 1053 Blindern, NO-0316 Oslo, Norway. {\texttt E-mail:
\{heidisu,ragnip\}@math.uio.no}}}

\maketitle

\begin{abstract} 
Polar varieties have in recent years been used by  Bank,
Giusti, Heintz, Mbakop, and Pardo, and by Safey El Din and Schost, to find
efficient procedures for determining points on all real components
of a given
non-singular algebraic variety. In this note we review the classical notion of
polars and polar varieties,
as well as the
construction of what we here call reciprocal polar varieties. In particular we
consider the case of real affine plane curves, and we give conditions for when
the polar varieties of singular curves contain points on all real components.
\end{abstract}


{\AMSMOS Primary 14H50; Secondary 14J70,14P05, 14Q05.
\endAMSMOS}

\section{Introduction}

The general study of polar varieties goes back to Severi, though polars of
plane
curves were introduced already by Poncelet. Polar varieties
of complex, projective, non-singular varieties have later been studied by many,
and the theory has been extended to singular varieties (see
\cite{P} and references therein).
We shall refer to these polar varieties as \emph{classical}
polar varieties.

Another kind of polar varieties, which we here will call \emph{reciprocal} polar
varieties, were introduced in \cite{bank04} under the name of dual
polar varieties. The definition involves a quadric hypersurface and polarity
of
 linear spaces with respect to this quadric. Classically, the
\emph{reciprocal curve} of a plane curve was defined to be the curve
consisting of the polar points of the 
tangent lines of the curve with respect to a given conic. 
The reciprocal curve is isomorphic to  the dual curve in $(\mathbb P^2)^\vee$, via the
isomorphism of $\mathbb P^2$ and $(\mathbb P^2)^\vee$ given by the
quadratic form defining the conic. 
 
Bank, Giusti, Heintz, Mbakop, and Pardo have proved
\cite{bank97,bank01,bank04,bank05} that polar varieties of
\emph{real, affine} non-singular varieties (with some requirements) contain
points on each connected
component of the variety, and  this property is useful in CAGD for finding a
point on each component. Related work has been done by Safey El Din and Schost
\cite{safey, schost}.

We will in this paper determine in which cases the polar varieties of a real
affine \emph{singular} plane curve contain at least one non-singular point of each
component of the curve.    In the next section we  briefly review the
definitions
of classical and reciprocal polar varieties, and state and partly prove 
versions of some results from \cite{bank04} adapted to our situation. The third
section treats the case of plane curves; we show that the presence of ordinary
multiple points does not affect these results, but that the presence of arbitrary
singularities does.

\section{Polar varieties}

Let $V\subset \mathbb{P}^n$ be a complex projective variety. Given a hyperplane
$H$ we can consider the affine space $\mathbb{A}^n:=\mathbb{P}^n \setminus H$,
where $H$ is called the hyperplane at infinity. We define the
corresponding affine variety $S$ to be the variety $V\cap \mathbb{A}^n$. We let
$V_{\mathbb{R}}$ and $S_{\mathbb{R}}$ denote the corresponding real varieties.  

If $L_1$ and $L_2$ are linear varieties in a projective space $\mathbb{P}^n$,
we
let  $\langle L_1,L_2\rangle$ denote the linear variety spanned by them. We say
that $L_1$ and $L_2$ intersect transversally if  $\langle
L_1,L_2\rangle=\mathbb{P}^n$; if they do not intersect transversally,
we write $L_1\not\pitchfork L_2$. 

If $I_{1}$ and $I_{2}$ are sub-vector spaces of $\mathbb A^n$ (considered as a
vector space), we say
that $I_1$ and $I_2$ intersect transversally if  $I_1+I_2=\mathbb{A}^n$;
if they do not intersect transversally, we write $I_1\not\pitchfork I_2$. 

\subsection{Classical polar varieties} 

Let us recall the definition of the classical polar varieties (or loci) of a
possibly singular variety. Consider
a flag of linear varieties in $\mathbb{P}^n$,
\begin{displaymath}
\mathcal{L}: L_0\subset L_1 \subset \ldots \subset L_{n-1}\subset \mathbb{P}^n
\end{displaymath} 
and a (non degenerate) variety $V\subset \mathbb{P}^n$ of codimension $p$.
Let $V_{\ns}$ denote the set of nonsingular points of $V$. For each point
$P\in V_{\ns}$, let $T_{P}V$ denote the \emph{projective} tangent space to $V$
at $P$.
The $i$-th polar variety ${W}_{L_{i+p-2}}(V)$,
$1\leq i \leq n-p$, of $V$ with respect to $\mathcal L$ is the Zariski closure of the
set
\begin{displaymath}
\{P \in V_{\ns}\setminus L_{i+p-2}\, |\, T_P V \not\pitchfork L_{i+p-2}\}.
\end{displaymath}

Take $H=L_{n-1}$ to be the hyperplane at infinity and let $S=V\cap
\mathbb{A}^n$ be the affine part of $V$. Then we can define the \emph{affine}
polar
varieties of $S$ with respect to the flag $\mathcal{L}$ as follows: 
 the $i$-th affine
polar variety
${W}_{L_{i+p-2}}(S)$ of $S$ with respect to $\mathcal L$ is the intersection
${W}_{L_{i+p-2}}(V)\cap
\mathbb{A}^n$.

Since $L_{n-1}$ is the hyperplane at infinity, and all the other elements of the
flag $\mathcal{L}$ are contained in $L_{n-1}$, we can look at the affine cone
over each element $L_j$ of the flag, considered  as a $(j+1)$-dimensional linear
sub-vector space $I_{j+1}$ of $\mathbb{A}^n$. Hence we get the flag (of vector
spaces)
\begin{displaymath}
\mathcal{I}: I_1\subset I_2\subset \ldots \subset I_{n-1} \subset \mathbb{A}^n,
\end{displaymath}
 and the affine polar variety ${W}_{L_{j-1}}(S)$ can also be defined as the
closure of the set
\begin{displaymath}
\{ P\in S_{\ns}\,|\, t_P S \not\pitchfork I_j\},
\end{displaymath} 
where $t_{P}S$ denotes the  \emph{affine} tangent
space to $S$ at $P$ translated to the origin (hence considered as a sub-vector
space of $\mathbb A^n$).
For more details on these two equivalent definitions of affine polar varieties,
see \cite{bank04}.  
\medskip
\begin{figure}[htbp!]
\caption{A conic with two tangents. The polar locus of the conic (with respect to
the point
of intersection of the two tangents) consists of the two points of tangency. The
\emph{polar} of the point is the line through the two points of tangency.}
  \label{polarconic}
  \begin{center}
    \includegraphics
   [angle=0,width=0.4\textwidth]
    {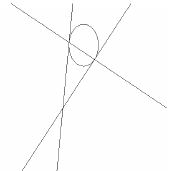}
  \end{center}
\end{figure}

If $f_{1},\ldots ,f_{r}$ are homogeneous polynomials in $n+1$ variables, we
let $\mathcal V(f_{1},\ldots,f_{r})\subset \mathbb P^n$ denote the corresponding
algebraic
variety (their common zero set). Similarly, if $f_{1},\ldots ,f_{r}$ are
polynomials in $n$ variables, we denote by $\mathcal
V(f_{1},\ldots,f_{r})\subset
\mathbb A^n$ the corresponding affine variety.
We shall often use the flag
\begin{equation}\label{flag}
\mathcal{L}: \mathcal{V}(X_0,X_2,\ldots,X_n)\subset
\ldots \subset
\mathcal{V}(X_0,X_n)\subset \mathcal{V}(X_0)\subset \mathbb P^n
\end{equation} 
and the corresponding affine flag (considered as vector spaces)
\begin{equation}\label{aflag}
\mathcal{I}: \mathcal{V}(X_2, \ldots, X_n)\subset 
\ldots \subset \mathcal{V}(X_{n-1},X_n)\subset \mathcal{V}(X_n)\subset
\mathbb A^n.
\end{equation}

The following proposition is part of a proposition stated and proved
in
\cite[2.4, p. 134]{bank01}. Since this particular result is valid under
fewer
assumptions, we shall give a different and more intuitive
proof.  We let $L$
denote any member of the flag $\mathcal{L}$. 
\medskip

\begin{proposition}\label{prop1} Let $S\subset \mathbb{A}^n$ be a  pure
$p$-codimensional reduced variety. Suppose that $S_{\mathbb{R}}$ is non-empty,
pure
$p$-codimensional, non-singular, and compact. Then
${W}_{L}(S_{\mathbb{R}})$ contains at least one point of each connected
component of $S_{\mathbb{R}}$.
\end{proposition}

\begin{proof} By an affine coordinate change, we may assume that the flag
$\mathcal{I}$ is the flag (\ref{aflag}).
The polar varieties form a sequence of inclusions
\begin{displaymath} {W}_{L_{n-2}}(S_{\mathbb{R}}) \subset
{W}_{L_{n-3}}(S_{\mathbb{R}}) \subset \ldots \subset
{W}_{L_{p-1}}(S_{\mathbb{R}}),
\end{displaymath} 
so it is sufficient to find a point of
${W}_{L_{n-2}}(S_{\mathbb{R}})$ on each component of the variety
$S_{\mathbb{R}}$. We will show that the maximum point for the last coordinate
$X_n$ of each component of $S_{\mathbb{R}}$ is also a point on
${W}_{L_{n-2}}(S_{\mathbb{R}})$.

The variety ${W}_{L_{n-2}}(S_{\mathbb{R}})$ is the set
\begin{displaymath}
\{P \in S_{\mathbb{R}}\, |\, t_P S \not\pitchfork \mathcal{V}(X_n)\}.
\end{displaymath}
Let $C$ be any connected component of $S_{\mathbb{R}}$, and let $A:=(a_1,\ldots,
a_n) \in C$ be a local maximum point for the last coordinate $X_n$ in  $C$. Such
a maximum point exists since $C$ is compact. Since $A$ is a local maximum point
for the $X_n$-coordinate, the variety $S_{\mathbb{R}}$ has to flatten out in the
$X_n$-direction in the point $A$, which means that the tangent space $t_A S$ is
contained in the hyperplane $\mathcal{V}(X_n)$. To show this, consider a real
local parameterization $S_{\mathbb{R}}$ in a neighborhood of $A$. The
neighborhood is chosen to be so small that it contains no other local maximum
points for the $X_n$-coordinate. 

Consider a local parameterization
\begin{displaymath}
\begin{array}{ccc} (s_1,\ldots,s_{n-p}) & \longmapsto &
(X_1(s_1,\ldots,s_{n-p}), \ldots, X_n(s_1,\ldots,s_{n-p})),
\end{array}
\end{displaymath} 
with
\begin{displaymath} A=(X_1(0,\ldots,0),\ldots,X_n(0,\ldots,0)).
\end{displaymath} 
The rows of the matrix
\begin{displaymath}
\left[
\begin{array}{ccc}
\frac{
\partial X_1}{
\partial s_1}(0,\ldots,0) & \cdots & \frac{
\partial X_n}{
\partial s_1}(0,\ldots,0) \\
\vdots & \ddots & \vdots \\
\frac{
\partial X_1}{
\partial s_{n-p}}(0,\ldots,0) & \cdots & \frac{
\partial X_n}{\partial s_{n-p}}(0,\ldots,0)
       \end{array}
\right]
\end{displaymath} 
span the tangent space $t_{A}S$ of $S_{\mathbb{R}}$ at $A$. We
will show that $\frac{
\partial X_n}{
\partial s_i}(0,\ldots,0)=0$ for
$i=1,\ldots,n-p$. Then we know that $t_A S$ is contained in
the hyperplane  $\mathcal{V}(X_n)$.

By the definition of derivatives we have that
\begin{displaymath}
\begin{array}{ccc}
\frac{
\partial X_n}{
\partial s_i}(0,\ldots,0)&=& \lim_{s_i \rightarrow
0^{-}}\frac{X_n(0,\ldots,0,s_i,0,\ldots,0)-X_n(0,\ldots,0)}{s_i-0}\\
\\
&=& \lim_{s_i \rightarrow
0^{+}}\frac{X_n(0,\ldots,0,s_i,0,\ldots,0)-X_n(0,\ldots,0)}{s_i-0}
\end{array}
\end{displaymath} 
If $s_i$ goes to $0$ from below, then $s_i$ is negative, and
$X_n(0,\ldots,0,s_i,0,\ldots,0)-X_n(0,\ldots,0)\leq 0$ since $X_n(0,\ldots,0)$
is the local maximum, so we have
\begin{displaymath}
\textstyle
\lim_{s_i \rightarrow
0^{-}}\frac{X_n(0,\ldots,0,s_i,0,\ldots,0)-X_n(0,\ldots,0)}{s_i-0} \geq 0
\end{displaymath} 
When $s_i$ goes to $0$ from above, then $s_i$ is positive, so
\begin{displaymath}
\textstyle
\lim_{s_i \rightarrow
0^{+}}\frac{X_n(0,\ldots,0,s_i,0,\ldots,0)-X_n(0,\ldots,0)}{s_i-0} \leq 0
\end{displaymath} 
Since these limits are equal, they both have to be zero, so
$\frac{\partial X_n}{\partial s_i}(0,\ldots,0)=0$ for $i=1,\ldots,n-p$.

Since $t_A S$ is contained in $\mathcal{V}(X_n)$, we have  $t_A S \not\pitchfork
\mathcal{V}(X_n)$, and hence $A$ is a point of ${W}_{L_{n-2}}(S_{\mathbb{R}})$.
\end{proof} 
\medskip

In the proof of this proposition we can replace the phrase ``local
maximum point of the $X_n$-coordinate" with ``local minimum point \ldots ",
since
the consideration we did for the tangent space is the same in both cases. So we
have found two points of each connected component, unless  the local maximum and
the local minimum for the last coordinate coincide. The local maximum and the
local minimum  coincide if and only if the variety is contained in the hyperplane
$\mathcal{V}(X_n)$, which implies that $V$ is degenerate.   

In the above proposition we did not assume that the coordinates were in generic
position with respect to the polynomials generating the variety, so we allow
situations where the polar variety can contain a piece of a component of
dimension greater than zero.

Safey El Din and Schost \cite{schost} show that one can find  points on each
component of a smooth, not necessarily compact, affine real variety. They use
all the polar varieties given by a flag, and intersect each of these varieties
with other varieties. The union of these intersections will be zero-dimensional
and it contains points from each connected component of the variety we started
with. 

\subsection{Reciprocal polar varieties}

Let $Q=\mathcal{V}(q)$ be a non-degenerate hyperquadric defined in $\mathbb{P}^n$ by a
polynomial $q$. If $A$ is a point, its polar
hyperplane $A^\perp$ with respect to $Q$ is the linear span of the points on
$Q$ such that the tangent hyperplanes to $Q$ at these points pass through $A$.
This gives the hyperplane $\sum \frac{\partial q}{\partial X_{i}}(A) X_{i}=0$.

If $H$ is a hyperplane, then its polar point, $H^\perp$, is the
intersection of the tangent hyperplanes for $Q$ at the points on $Q\cap H$.
Finally, if $L$ is a linear space of dimension $d$, its polar space
$L^\perp$ is the intersection of the polar hyperplanes to points in $L$.
 Equivalently, $L^{\perp}$ is the linear span of all points $H^{\perp}$,
where $H$ is a hyperplane containing $L$. The dimension of $L^{\perp}$ is
$n-d-1$.
\medskip

As an example, consider the case $n=3$. If $A$ is a point in $\mathbb{P}^3$
then $A^{\perp}$ is the plane in $\mathbb{P}^3$ defined by the polynomial
$\frac{\partial q}{\partial X_0}(A)X_0 + \ldots + \frac{\partial q}{\partial
X_3}(A)X_3$. If $L$ is a hyperplane defined by the polynomial
$b_0X_0+b_1X_1+b_2X_2+b_3X_3$, then $L^{\perp}$ is the point $A$ such that 
\begin{displaymath} \textstyle
(b_0:b_1:b_2:b_3)=(\frac{
\partial q}{
\partial X_0}(A): \frac{
\partial q}{
\partial X_1}(A):\frac{
\partial q}{
\partial X_2}(A):\frac{
\partial q}{
\partial X_3}(A))
\end{displaymath} 
Finally, if $L$ is the line spanned by the points $A$ and $B$,
then $L^{\perp}$ is the intersection of the two hyperplanes $A^{\perp}$ and
$B^{\perp}$, i. e.,
\begin{displaymath} \textstyle
 L^{\perp}=\mathcal{V}(\sum_{i=0}^{3}\frac{
\partial q}{
\partial X_i}(A)X_i)\cap\mathcal{V}(\sum_{i=0}^{3}\frac{
\partial q}{
\partial X_i}(B)X_i) 
\end{displaymath}

Note that if $Q$ is defined by the polynomial $q=\sum_{i=0}^{n}X_i^2$, then
the polar variety of a point $(a_0: \ldots :a_n)\in {\mathbb{P}^n}$ is
the hyperplane 
$H=\mathcal{V}(\sum_{i=0}^{n}a_iX_i) \subset {\mathbb{P}^n}$.
\medskip

Let $\mathcal{L}: L_0 \subset L_1 \subset \dots \subset L_{n-1}$ be a flag in
$\mathbb{P}^n$, where $L_{n-1}$ is the hyperplane $H$ at infinity if we consider the
affine space. We then get the polar flag with respect to $Q$:
\begin{displaymath}
\mathcal{L^{\perp}}: L_{n-1}^{\perp}\subset L_{n-2}^{\perp}\subset \dots \subset
L_1^{\perp}\subset L_0^{\perp}
\end{displaymath} 

\begin{definition}[cf. {\cite[p. 527]{bank04}}] The i-th reciprocal polar variety
 $W_{L_{i+p-1}}^{\perp}(V)$, $1\leq i \leq n-p$, of
a variety $V$ with respect to the flag $\mathcal{L}$, 
is defined to be the Zariski
closure of the set
\begin{displaymath}
\{ P\in V_{\ns}\setminus L_{i+p-1}^{\perp}\,|\, T_P V \not\pitchfork
\langle P,L_{i+p-1}^{\perp}\rangle ^{\perp}\}
\end{displaymath}
\end{definition}
\medskip

When $V$ is a hypersurface in $\mathbb{P}^n$,  the reciprocal
polar
variety $W_{L_{n-1}}^{\perp}(V)$ is the set $\clos\{ P\in V_{\ns}\setminus
L_{n-1}^{\perp}\,|\, T_P V \supset \langle P,L_{n-1}^{\perp}\rangle ^{\perp}
\}$. In this case, $\langle
P,L_{n-1}^{\perp}\rangle$ is the line spanned by $P$ and the point $L_{n-1}^{\perp}$.
Since $A^{\perp \perp}=A$, and $A \subseteq B$ implies $A^\perp \supseteq
B^\perp$, it follows that 
\begin{displaymath} 
T_P V \supseteq \langle P,L_{n-1}^{\perp}\rangle ^{\perp}
\Leftrightarrow {T_P V}^\perp \in \langle P,L_{n-1}^{\perp}\rangle.
\end{displaymath} 
The point ${T_P V}^\perp$ is on
the line $\langle P,L_{n-1}^{\perp}\rangle$ if and only if the point $L_{n-1}^{\perp}$
is on the line $\langle P,{T_P V}^\perp\rangle$. So, when $V$ is a hypersurface,
the $(n-1)$-th reciprocal variety is
\begin{displaymath} W_{L_{n-1}}^{\perp}(V)=\clos\{ P\in V_{\ns}\setminus
(\{L_{n-1}^{\perp}\}\cup H)\,|\,L_{n-1}^{\perp} \in \langle P,{T_P
V}^\perp\rangle\}.
\end{displaymath}
This way of writing the reciprocal polar variety can sometimes be
useful, and it gives at better geometric understanding of the reciprocal polar
variety.
\medskip

Let $S \subset \mathbb{A}^n=\mathbb P^n\setminus H$ denote the affine part of the
variety $V$, 
where $H=L_{n-1}$ is the hyperplane at infinity, and let $L\in \mathcal L$.
We define the \emph{affine} reciprocal polar variety 
$W_L^{\perp}(S)$ to be the affine part
$W_L^{\perp}(V)\cap \mathbb{A}^n$ of $W_L^{\perp}(V)$. 
The linear variety $\langle
P,L^{\perp}\rangle ^{\perp}$ is contained in the hyperplane at infinity, so we
can consider the affine cone of $\langle P,L^{\perp}\rangle ^{\perp}$ as a
linear variety $I_{P,L^{\perp}}$ in the affine space. Then the affine reciprocal
polar variety can be written as
\begin{displaymath} W_L^{\perp}(S)=\clos \{P \in S_{\ns}\setminus
L^{\perp}\,| \, t_P S \not\pitchfork I_{P,L^{\perp}} \}
\end{displaymath} where $t_P S$ is the affine tangent space at $P$, translated to the
origin.
\medskip

The following result, formulated slightly differently, is proved
in \cite{bank04}.
\medskip

\begin{proposition}[{\cite[p. 529]{bank04}}] \label{thm2}  
Let $S_{\mathbb{R}}\subset \mathbb A_{\mathbb R}^n$ be
a non-empty, non-singular real variety of pure codimension
$p$. Let $\mathcal{L}$ be a flag in $\mathbb P^n$, where $L_{n-1}=H=
\mathbb P^n\setminus \mathbb A^n$ is the
hyperplane at infinity. Assume
$Q=\mathcal V(q)$, where $q$ restricts to a positive definite quadratic form on 
$\mathbb A_{\mathbb R}^n$. Assume $H^{\perp}\notin S_{\mathbb{R}}$.
 Let $L$ be any member of the flag
$\mathcal L$ with $\dim L \ge p$. Then the real
affine reciprocal polar variety $W_L^{\perp}(S_{\mathbb{R}})$ 
contains at least one point from each connected component of $S_{\mathbb{R}}$. 
\end{proposition} 
\medskip

In Proposition \ref{thm2} we had to
assume that $H^{\perp}\notin S_{\mathbb{R}}$ in order
to prove that $W_H^{\perp}(S_{\mathbb{R}})$ contains a point from each component
of $S_{\mathbb{R}}$.
The following proposition
states
that if the variety is a hypersurface and contains $H^{\perp}$, then we can choose
another quadric $Q^{\prime}$, so that the polar point $H^{\perp ^{\prime}}$
with respect to this quadric is not on $S_{\mathbb{R}}$, and we will thus still
be able to find points on each component of the hypersurface.
\medskip

\begin{proposition}\label{transform} Let $V\subset \mathbb P^n$ be a
hypersurface and $H$ a hyperplane, and set $\mathbb A^n=\mathbb
P^n\setminus H$ and
$S= V\cap \mathbb A^n$.
Asssume
$S_{\mathbb{R}}$ is
non-empty and non-singular.
Given any 
point $A\in \mathbb{A}_{\mathbb{R}}^n \setminus S_{\mathbb{R}}$, there exists a
quadric $Q^{\prime}$ such that the polar point $H^{\perp^{\prime}}$ with respect to
$Q^{\prime}$ is equal to $A$, and such that
the real affine reciprocal polar variety
$W_H^{\perp^{\prime}}(S_{\mathbb{R}})$  contains at least one point from each
connected component of $S_\mathbb{R}$.
\end{proposition}

\begin{proof}
By Proposition \ref{thm2} it suffices to show that we can find $Q^{\prime}=\mathcal V(q^{\prime})$ 
such that $H^{\perp^{\prime}}=A=(1:a_{1}:\ldots :a_{n})$ and such that $q^{\prime}$ restricts to 
a positive definite form on $\mathbb A^n$. Take $q^{\prime}=(1+\sum_{i=1}^na_{i}^2)X_{0}^2
-2\sum_{i=1}^na_{i}X_{0}X_{i}+\sum_{i=1}^nX_{i}^2$. 
Then $\frac{\partial q^{\prime}}{\partial X_{0}}(A) =2$, and 
$\frac{\partial q^{\prime}}{\partial X_{i}}(A) =0$ for $i=1,\ldots,n$. The restriction
of $q^{\prime}$ to $\mathbb A^n$ is $\sum_{i=1}^nX_{i}^2$.
\end{proof}
\medskip

The next proposition provides an explicit description of the affine 
reciprocal polar variety of a hypersurface $S\subset \mathbb{A}^n$. 
A similar result is proved in the more general case of a complete intersection
variety in \cite[3.1, p. 531]{bank04}.
\medskip

\begin{proposition}\label{minors}
Let $V=\mathcal{V}(f)\subset \mathbb P^n$ be a real
hypersurface, and let $Q=\mathcal{V}(q)$ be defined by an irreducible
quadratic polynomial $q$.
We consider the affine space $\mathbb A^n=\mathbb P^n\setminus H$, where
$H=\mathcal{V}(X_0)$.
Then the
affine reciprocal polar variety  $W_{H}^{\perp}(S)$ of $S=V\cap\mathbb{A}^n$
with respect to $Q$ is equal to the closure of the  intersection of 
$S_{\ns}\setminus \{H^{\perp}\}$
with the variety defined by the $2$-minors of the matrix
\begin{displaymath}
\left( \begin{array}{ccc} \frac{\partial f}{\partial X_1} & \cdots & 
\frac{\partial f}{\partial X_n} \\
\frac{\partial q}{\partial X_1}&\cdots &  \frac{\partial q}{\partial X_n}
\end{array}
\right).
\end{displaymath}
\end{proposition}

\begin{proof} The reciprocal variety $W_{H}^{\perp}(S)$ is the closure of the set
\begin{displaymath}
\{ P\in S_{\ns}\setminus H^{\perp}\,|\,I_{P,H^{\perp}}\not\pitchfork t_{P}S\}.
\end{displaymath} 
We want to  show that for $P=(1:p_1:\ldots:p_n)\in S$ the affine cone 
$I_{P,H^{\perp}}$  over $\langle P,H^{\perp} \rangle
^{\perp}$ intersects the tangent space $t_P S$ non-transversally if and only if
\begin{displaymath}
\rank \left( \begin{array}{ccc} \frac{\partial f}{\partial X_1}(P) & \cdots &  \frac{\partial f}{\partial X_n}(P) \\
 \frac{\partial q}{\partial X_1}(P)&\cdots &  \frac{\partial q}{\partial X_n}(P)
\end{array}
\right)\leq 1.
\end{displaymath}
Note that we have 
\[
\langle P,H^{\perp}\rangle ^{\perp}=P^{\perp}\cap H,
\]
so, since $P^{\perp}=\mathcal V(\sum_{i=0}^{n}\frac{\partial q}{\partial
X_i}(P)X_{i})$ and $H=\mathcal V(X_{0})$, we find
\begin{displaymath}
\langle P,H^{\perp}\rangle ^{\perp}=\{(0:X_1:\ldots:X_n)\,|\,\sum_{i=1}^{n}\frac{\partial q}{\partial X_i}(P)X_{i}=0\}
\end{displaymath}
The affine cone $I_{P,H^{\perp}}$ is the hyperplane
$\mathcal{V}(\sum_{i=1}^{n}\frac{\partial q}{\partial X_i}(P)X_{i})$. 
The affine
tangent space $t_P S$ is the hyperplane given by $
\sum_{i=1}^n\frac{\partial f}{\partial X_i}(P)X_{i}=0$, which implies what we
wanted to prove.
\end{proof}
\medskip

Note that 
the square of the distance function in the affine space
$\mathbb{A}^n$ is given
by a quadratic polynomial, so if we let $Q=\mathcal V(q)$ be a quadric such 
that $q$ restricts to this
polynomial, we see that the variety defined in  Theorem 2 in \cite{safey} is
nothing but the affine reciprocal polar variety with respect to $Q$. (This is
because elimination of the extra variable gives the $2$-minors of the matrix
in Proposition \ref{minors}.) Hence \cite[Thm. 2]{safey} follows from
Proposition \ref{thm2}.
\medskip

\section{Polar varieties of real singular curves}

In this section our varieties will be curves in $\mathbb{P}^2$, so the flags
will in this case be of the following form
\begin{displaymath}
\mathcal{L}: L_0 \subset L_1 \subset \mathbb{P}^2
\textrm{ and }
\mathcal{L}^{\perp}: L_1^{\perp} \subset L_0^{\perp} \subset \mathbb{P}^2
\end{displaymath} where $L_0$ is a point, and $L_1$ is the line at infinity,
when
we look at the affine case. This gives only one interesting classical polar
variety and reciprocal polar variety for each choice of flag, namely
$W_{L_0}(V)$ and $W_{L_1}^{\perp}(V)$. 
\medskip

Let $V=\mathcal V(f)\subset \mathbb P^2$ be a plane curve.
The \emph{polar} of $V$  with respect
to a point $A=(a_{0}:a_{1}:a_{2})$
is the curve $V^{\prime}=\mathcal V(\sum_{i=0}^2a_{i}\frac{\partial f}{\partial
X_{i}})$.
It follows from the definition that the polar variety of $V$ with respect to $A$
is equal to 
\[
W_{A}(V)=V_{\ns}\cap V^{\prime}.\]
Given a conic $Q=\mathcal V(q)$ and a line $L$, the \emph{reciprocal polar} of
the affine curve $S=V\cap \mathbb A^2\subset \mathbb A^2=\mathbb P^2\setminus L$,
is the curve $S^{\prime \prime}=\mathcal V(\frac{
\partial f}{
\partial X_{1}} \frac{
\partial q}{
\partial X_{2}}-
\frac{\partial f}{\partial X_{2}} \frac{\partial q}{
\partial X_{1}})$, where $f$ and $q$ are dehomogenized by setting $X_{0}=1$.
The affine reciprocal polar variety of $S$ is equal to
\[W_{L}^{\perp}(S)=
S_{\ns}\cap S^{\prime \prime}.\]
More generally, we define the \emph{reciprocal polar} of the curve $V\subset
\mathbb P^2$ with respect to $Q$ and a point $A=(a_{0}:a_{1}:a_{2})$ to be the
curve 
$V^{\prime \prime}=\mathcal V(\det(f,q,A))$,
where $\det(f,q,A)$ is the determinant of the matrix
\[
\left[\begin{array}{ccc}
a_{0} & a_{1} & a_{2} \\
\frac{\partial f}{\partial X_0} & \frac{\partial f}{\partial X_1} & 
\frac{\partial f}{\partial X_2} \\
\frac{\partial q}{\partial X_0} & \frac{\partial q}{\partial X_1} & 
\frac{\partial q}{\partial X_2} 
       \end{array}
\right].
\]
The reciprocal polar of the affine curve $S$ is then the affine part of the
reciprocal polar with respect to the origin $(1:0:0)$.

Note that the reciprocal polars form a linear system on $\mathbb P^2$ of degree $d$,
whereas the classical polars form a linear system of degree $d-1$. A
classical polar variety consists of at most $d(d-1)$ points, whereas a
reciprocal polar variety can have as many as $d^2$ points.
\medskip

In the following sections we will consider affine curves with singularities, and
we determine in which cases the classical polar variety or the reciprocal polar
variety will contain non-singular points of each connected component of the
real part of the curve.

\subsection{Classical polar varieties of real singular affine curves} 
The classical
polar variety of an affine curve $S\subset \mathbb{A}^2$ associated to a given
flag $\mathcal L$ is the set
\begin{displaymath}
W_{L_{0}}(S)=\{ P\in S_{\ns}\,|\, I_1 = t_P S\},
\end{displaymath} 
where $I_{1}$ is the affine line equal to the cone over the point $L_{0}$.
By an \emph{ordinary real multiple point} we shall mean a singular point
with at least two real branches, such that the tangent
lines intersect pairwise transversally. If the curve
$S_{\mathbb R}$ only has ordinary
real multiple points as singularities, we have the following proposition.
\medskip

\begin{proposition}\label{curve} Suppose 
$S_{\mathbb{R}}$ is non-empty and
compact and has only
ordinary real multiple points as singularities. Then ${W}_{L_{0}}(S_{\mathbb{R}})$
contains at least one non-singular point of each connected component.
\end{proposition}

\begin{proof} We may assume that $I_{1}=\mathcal{V}(X_2)$. Let $C$
be a connected component of $S_{\mathbb{R}}$. The component $C$ has a local
maximum point for the coordinate $X_2$, since $C$ is compact. If this point is
a non-singular point, then we know by the proof of Proposition
\ref{prop1} that it is contained in the real polar variety. 
Assume on the contrary that the maximum point $P$ is a singular point, hence an
ordinary real multiple point.
If $P=(p_1,p_2)$ is a local maximum point for the $X_2$-coordinate then each of
the real branches through $P$ has $p_2$ as a local maximum for the
$X_2$-coordinate, so the line $\mathcal{V}(X_2-p_2)$ is a tangent line for each
branch. This means that the branches have a common tangent line, hence they do
not intersect transversally.
\end{proof}
\medskip

When it comes to singularities other than ordinary  real multiple points, we
cannot
say whether the singularity can be a maximum point for the $X_2$-coordinate.
But we know that the local maximum and the local minimum points are on the real
polar $S^{\prime}_{\mathbb R}$. One point on a component cannot be both a minimum
and a maximum
unless the component is a line, which is not compact, and therefore excluded. If
the
singularity is a local maximum point, we know that the minimum point is also on
the polar. So we can allow each real component of the curve to have one
additional singularity which is not an ordinary real multiple point, and the
above result will still remain valid.

 For
curves with arbitrary singularities, clearly we can have a
situation where the points with the maximum and the minimum values for the last
coordinate are both singular, with the result that the conclusion of Proposition
\ref{curve} will not be valid for the given choice of flag $I_{1}$. One could ask
whether it is possible to choose a different affine flag (different affine
coordinates) so that the results still holds.
\medskip

\begin{proposition} There exists an affine singular plane curve $S_{\mathbb R}$
such that for no
choice of flag $I_{1}$ does the polar variety contain a non-singular
point from each component.
\end{proposition}

\begin{proof} We prove this by giving an example of such a curve. Since cusps
disturb the continuity of the curvature of a curve, it is natural to look for
examples among curves with cusps. So we want to find a curve with at least two
components, where each component has cusps. One way to construct such a curve is
to consider two, not necessarily irreducible, affine curves $\mathcal{V}(f)$ and
$\mathcal{V}(g)$, and then look at the curve
$\mathcal{V}(h)$, where $h=f^2+\epsilon g^3$, which will have cusps at
the
points of intersection between $\mathcal{V}(f)$ and $\mathcal{V}(g)$. We let
$\mathcal{V}(f)$ be the union of two disjoint circles, and $\mathcal{V}(g)$ the
union of four lines, where two of the lines intersect one of the circles twice,
and the other two lines intersect the other circle twice. We can take
\begin {small}
\begin{displaymath} 
f=(X_1^2+X_2^2-1)((X_1-4)^2+(X_2-2)^2-1) 
\end{displaymath}
\end{small} 
and 
\begin{small}
\begin{displaymath}
\textstyle
g=(X_2-\frac{1}{2})(X_2+\frac{1}{2})(X_1-\frac{7}{2})(X_1-\frac{9}{2}). 
\end{displaymath} 
\end{small}
The curve $\mathcal{V}(h)=\mathcal{V}(f^2+\frac{1}{100}g^3)$
has then four compact components with two cusps on each component. 

\begin{figure}[htbp!]
  \caption{The curves $\mathcal{V}(f)$ and $\mathcal{V}(g)$.}
  \label{sirkleroglinjer}
  \begin{center}
  \includegraphics{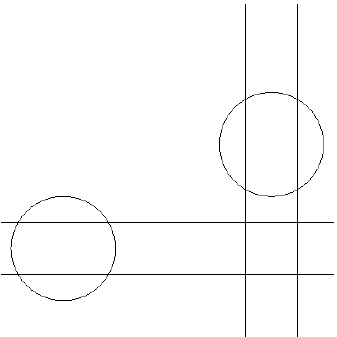}
  \end{center}
\end{figure}

\begin{figure}[htbp!]
  \caption{The curve $\mathcal{V}(f^2+\frac{1}{100}g^3)$.}
  \label{cuspkurve}
  \begin{center}
    \includegraphics{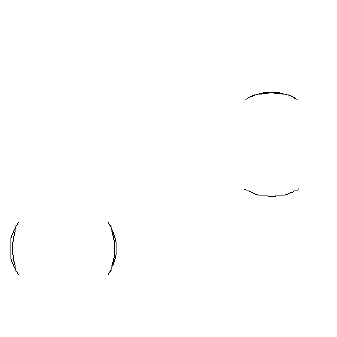}
  \end{center}
\end{figure}

A non-singular point $P$ is on the affine polar variety with respect to a line
 $I=\mathcal{V}(aX_1+bX_2)$ if the tangent line
$\mathcal{V}(\frac{\partial f}{\partial X_1}(P)X_1+\frac{\partial f}{
\partial X_2}(P)X_2)$ is equal to $I$.
Consider the
Gauss map $\gamma\colon S_{\mathbb R}\to \mathbb P^1_{\mathbb R}$ given by
\begin{displaymath} \textstyle
P\mapsto (\frac{
\partial f}{
\partial
  X_1}(P):\frac{
\partial f}{
\partial X_2}(P)).
\end{displaymath}
Let $C_i$, $i=1,\ldots,4$,
denote the connected components of $S_{\mathbb R}$. Let us show that 
$\cap_{i=1}^4\gamma(C_i)=\emptyset$.

The two lower components have cusps at the points
$(-\frac{\sqrt{3}}{2},\pm\frac{1}{2})$ and
$(\frac{\sqrt{3}}{2},\pm\frac{1}{2})$,
and the map $\gamma$ sends these points to the lines
$\mathcal{V}(3X_1\mp\sqrt{3}X_2)$ and
$\mathcal{V}(3X_1\pm\sqrt{3}X_2)$
respectively. Each of these components have points that are mapped to
the line $\mathcal{V}(X_1)$ by $\gamma$. The components do not
have any inflection points, so the image of the map $\gamma$ varies continuously
from $\mathcal{V}(X_1)$ to $\mathcal{V}(3X_1+\sqrt{3}X_2)$ and
$\mathcal{V}(3X_1-\sqrt{3}X_2)$, hence the other tangent lines lie in the sector
between the two lines $\mathcal{V}(3X_1+\sqrt{3}X_2)$ and
$\mathcal{V}(3X_1-\sqrt{3}X_2)$. 

The same happens if we calculate the tangent lines at the cusps for the two
upper components. The tangent lines of these components have to be in the sector
between  $\mathcal{V}(3X_2+\sqrt{3}X_1)$ and 
$\mathcal{V}(3X_2-\sqrt{3}X_1)$,
which contains the line $\mathcal{V}(X_2)$. These two sectors do not intersect
except at the origin, so for any choice of line $I=\mathcal{V}(aX_1+bX_2)$
the set $(\gamma ^{-1}(I))\cap C$ is empty or consists of points of the two
lower components or points of the two upper components.

\begin{figure}[htbp!]
  \caption{The lines $\mathcal{V}(3X_1+\sqrt{3}X_2)$,
$\mathcal{V}(3X_1-\sqrt{3}X_2)$, and $\mathcal{V}(X_1)$.}
  \label{linjermedV(x)}
  \begin{center}
    \includegraphics{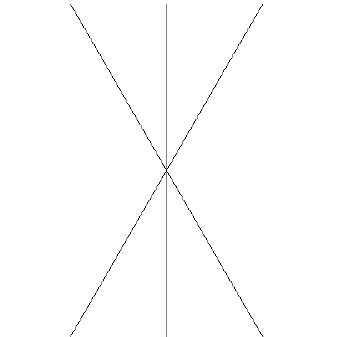}
  \end{center}
\end{figure}

\begin{figure}[htbp!]
  \caption{The lines $\mathcal{V}(3X_2+\sqrt{3}X_1)$,
$\mathcal{V}(3X_2-\sqrt{3}X_1)$, and $\mathcal{V}(X_2)$.}
  \label{linjermedV(y)}
  \begin{center}
    \includegraphics{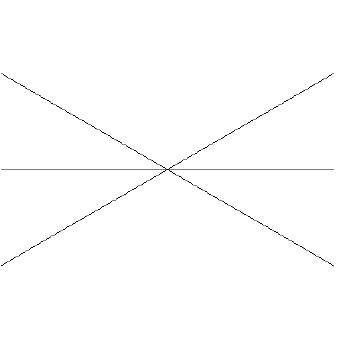}
  \end{center}
\end{figure}
\end{proof}
\medskip

\subsection{Reciprocal polar varieties of affine real singular curves}

As for the classical polar varieties we will here determine whether Proposition
\ref{thm2} is valid for affine curves with singularities. We
can easily prove the following for curves with ordinary multiple points.
\medskip

\begin{proposition}\label{proprpv} Let $S\subset \mathbb A^2=\mathbb
P^2\setminus L$
be an affine plane curve.
Suppose  that 
$S_{\mathbb{R}}$ is non-empty and has no other
singularities than ordinary real multiple points. Further, let $Q$ be defined
by a polynomial $q$ which restricts to a positive definite quadratic form on
$\mathbb A^2_{\mathbb R}$. Assume that $L^{\perp}$ 
is not contained in $S_{\mathbb{R}}$. Then the real
affine reciprocal polar variety $W_L^{\perp}(S_{\mathbb{R}})$ 
contains at least one non-singular point from each connected component of
$S_{\mathbb{R}}$. 
\end{proposition}

\begin{proof} Note that the hypotheses imply that $L^{\perp}\in \mathbb A^2=
\mathbb P^2\setminus L$. Moreover, the restriction of the quadratic polynomial
$q$ to $\mathbb A^2$ defines a distance function on $\mathbb A^2_{\mathbb R}$.
 The proof of Proposition
\ref{thm2} given in \cite{bank04}
consists, in this case,  of showing that for each component, a point with
the shortest distance to the point $L^{\perp}$  is  a point on the
reciprocal polar variety. So, we must show that if the component contains
ordinary real multiple points, then these points can not be
among the points with the shortest distance to
the  point $L^{\perp}$.

Assume on the contrary that, for a given component $C$ of the curve $S_{\mathbb
R}$, there is an ordinary real multiple point $P$ which has the
shortest distance to the point $L^{\perp}$.
The conic with centre in $L^{\perp}$ and radius dist$(L^{\perp},P)$ will by our
assumption be tangent to the component at the point $P$. 
Since $P$ is an ordinary real multiple point,
there is at least one other real branch, and this branch intersects the
conic
transversally. Hence this branch contains points inside the conic, and these
points are closer to $L^\perp$ than $P$.
\end{proof}
\medskip

As in the case of classical polar varieties we cannot prove the above
proposition for curves with arbitrary singularities, since the situation will
depend on the type of singularities and how the singularities are placed on the
component.
\medskip

\section{Examples}

In this section we shall look at some examples of singular plane affine curves
and
their polars and reciprocal polars, thus illustrating the propositions in the
previous sections. We use {\sc Surf} \cite{surf} to draw the curves. 
\medskip

{\sc Example 1.} Consider the real affine curve $S_{1}$ of degree $6$ 
defined by the polynomial
\begin{small}
\begin{displaymath}
\begin{array}{l}
f_{1}:=X_1^6+3X_1^4X_2^2-12X_1^4X_2+7X_1^4+3X_1^2X_2^4-24X_1^2X_2^3+66X_1^2X_2^2\\
-132X_1^2X_2+136X_1^2+X_2^6-12X_2^5+59X_2^4-132X_2^3+84X_2^2+144X_2-143.
\end{array}
\end{displaymath}
\end{small}
 This curve is compact and smooth, and it has three connected components.

\begin{figure}[htbp!]
\caption{The curve $S_{1}$.}
  \label{polareks1_1}
  \begin{center}
    \includegraphics{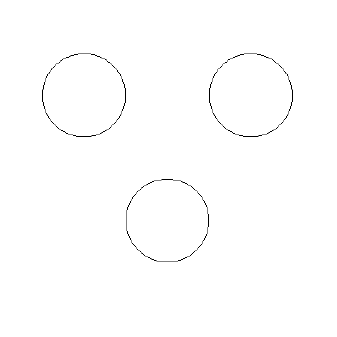}
  \end{center}
\end{figure}

We let $\mathcal{L}$ be the flag
\begin{displaymath}
L_{0}=\mathcal{V}(X_0,X_1)\subset L_{1}=\mathcal{V}(X_0).
\end{displaymath} The real polar variety is  the set
 \begin{displaymath}\textstyle
W_{L_{0}}((S_{1})_{\mathbb R})=\{P \in (S_{1})_{\mathbb{R}}\,|\,\frac{
\partial f_{1}}{
\partial X_2}(P)=0\},
\end{displaymath} 
which is equal to intersection of
$(S_{1})_{\mathbb{R}}$ and its polar $(S_{1}^{\prime})_{\mathbb R}=
\mathcal{V}(\frac{\partial
f_{1}}{
\partial X_2})_{\mathbb R}$. We see that the polar variety contains
points from
each connected component and that the points of the affine polar variety are
exactly those points on each component which give maximal and minimal values for
the $X_1$-coordinate. 

\begin{figure}[htbp!]
\caption{The curve $S_{1}$ and its polar $S_{1}^\prime$.}
  \label{polareks1_2}
  \begin{center}
   \includegraphics{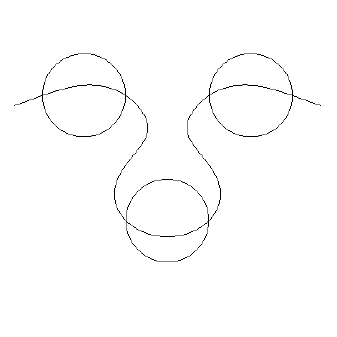}
  \end{center}
\end{figure}

Let $Q$ be the standard quadric, given by $q=\sum_{i=0}^2X_{i}=0$. The real
affine reciprocal polar variety consists of the real points of the
intersection between $S_{1}$ and its reciprocal polar $S_{1}^{\prime \prime}
=\mathcal V(X_2\frac{
\partial
  f_{1}}{
\partial X_1}-X_1\frac{
\partial f_{1}}{
\partial X_2})$, and we know from Proposition \ref{thm2} that also the real
affine reciprocal polar variety contains points from each connected component.
The reciprocal polar variety consists of the points on each component with the
locally shortest or longest Euclidean distance to the origin.

\begin{figure}[htbp!]
\caption{The curve $S_{1}$ and its reciprocal polar $S_{1}^{\prime \prime}$.}
  \label{polareks1_3}
  \begin{center}
    \includegraphics{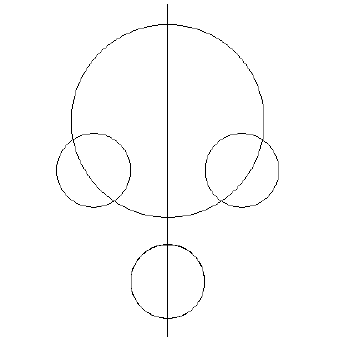}
  \end{center}
\end{figure}
\medskip

{\sc Example 2.} To illustrate Proposition \ref{curve} we consider an affine
irreducible curve $S_{2}$ with two compact components and with one ordinary
double point on each of the components. The curve is given by the polynomial
\begin{small}
\begin{displaymath}\textstyle
f_{2}=((X_1+2)X_2-(X_1+2)^6-X_2^6)(X_1X_2-X_1^6-X_2^6)+\frac{1}{100}X_2^6.  
\end{displaymath}
\end{small}

\begin{figure}[htbp!]
\caption{The curve $S_{2}$.}
  \label{dobbeltpktupolar}
  \begin{center}
   \includegraphics[angle=0,width=0.5\textwidth]{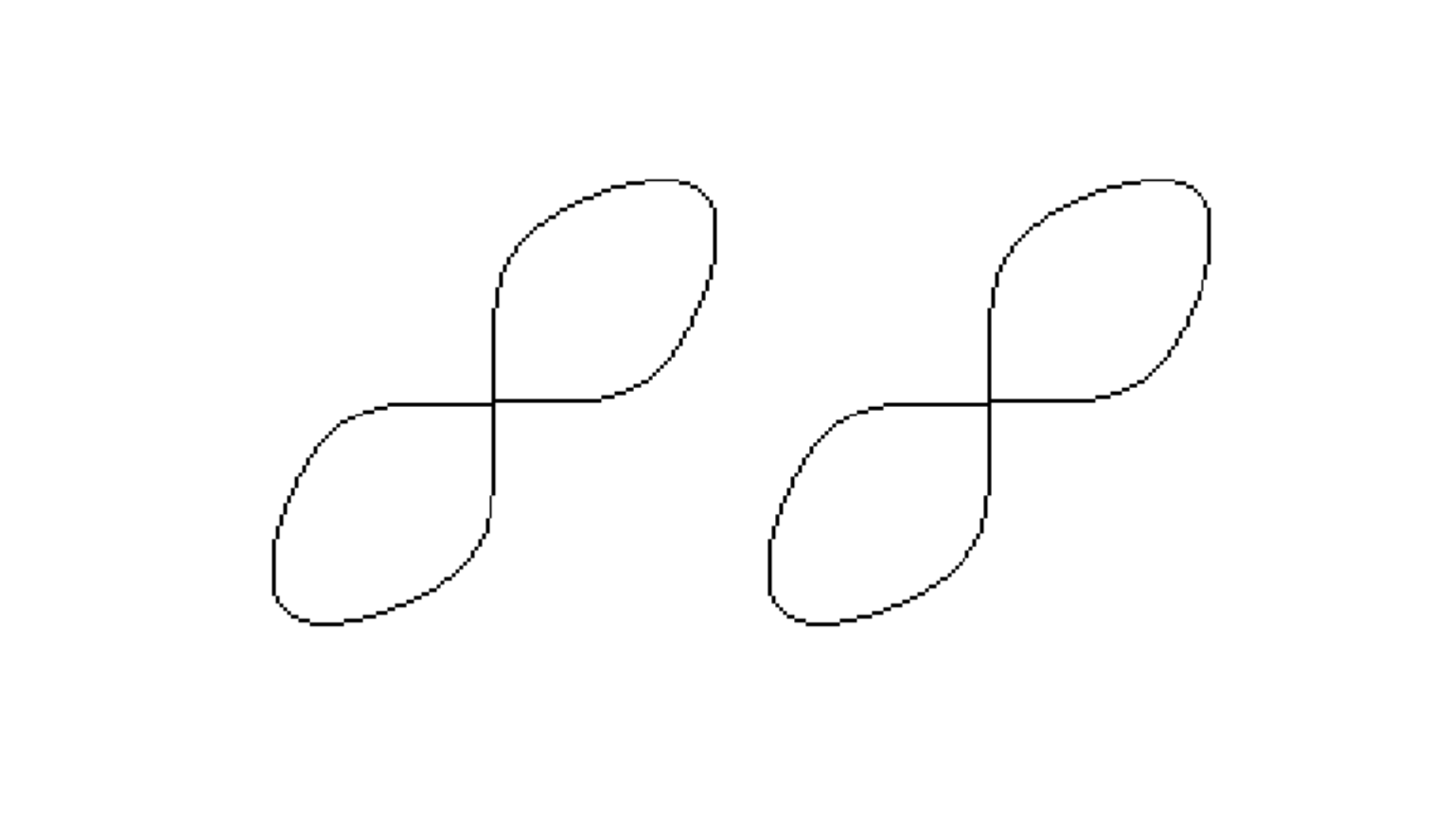}
  \end{center}
\end{figure}
We let $\mathcal{L}$ be as in the example above, and we see that
the affine polar variety contains non-singular points from each of the
components of the curve.

\begin{figure}[htbp!]
\caption{The curve $S_{2}$ and its polar $S_{2}^{\prime}$.}
  \label{dobbeltpktmpolar}
  \begin{center}
    \includegraphics[angle=0,width=0.5\textwidth]{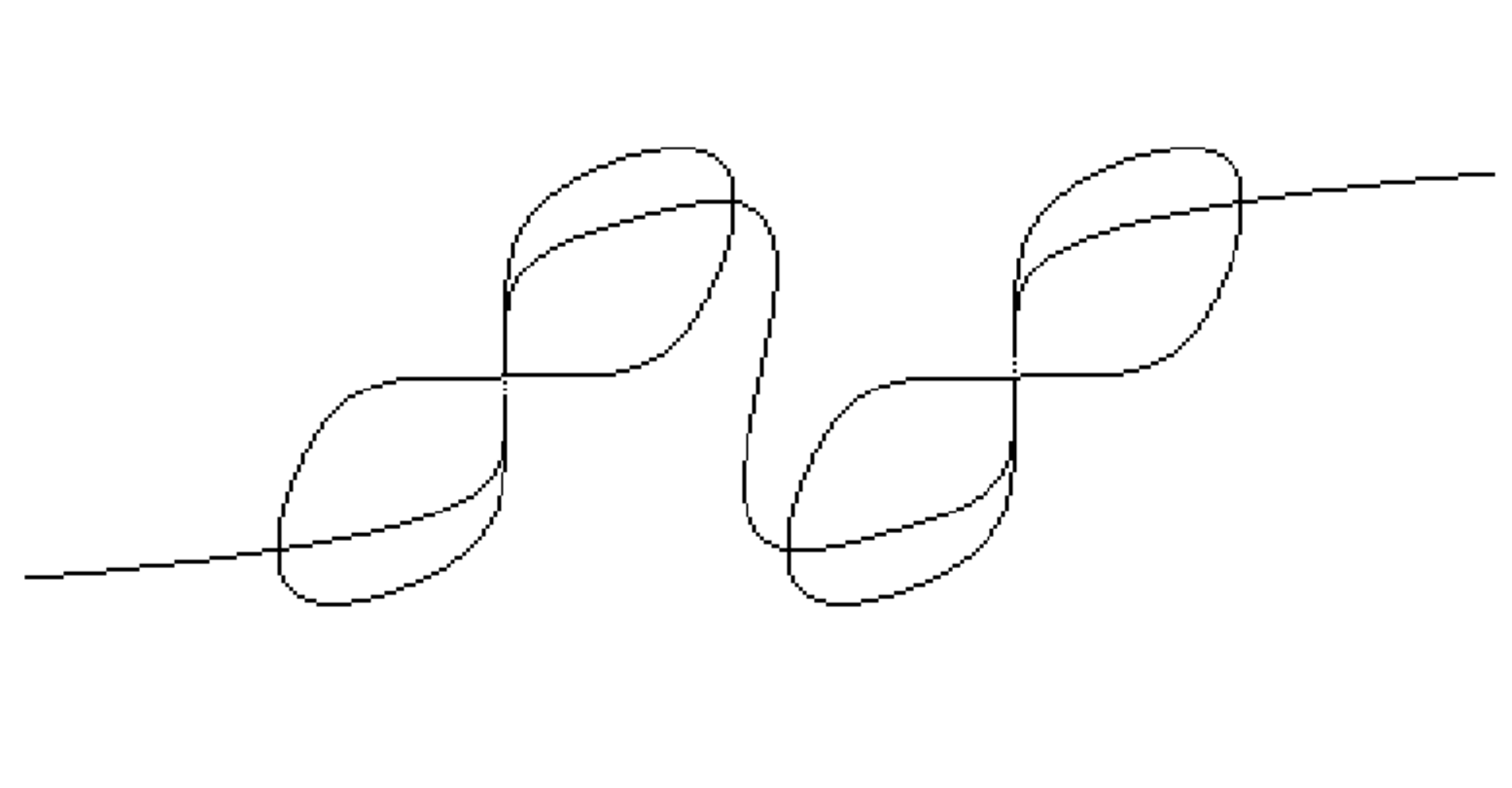}
  \end{center}
\end{figure}
\medskip

{\sc Example 3.} Consider the irreducible curve $S_{3}=\mathcal{V}(f_{3})$,
where
\begin{small}
\[
\textstyle 
f_{3}
=\mathcal{V}(144-24X_2^2-88X_1^2+X_2^4-X_1^6+17X_1^4-14X_2^2X_1^2%
+\frac{1}{100}X_2^6).\]
\end{small}
The real part of the curve consists of two non-compact components with one
ordinary double point on each component.

\begin{figure}[htbp!]
\caption{The curve $S_{3}$.}
  \label{ikkekomp2dp}
  \begin{center}
   \includegraphics{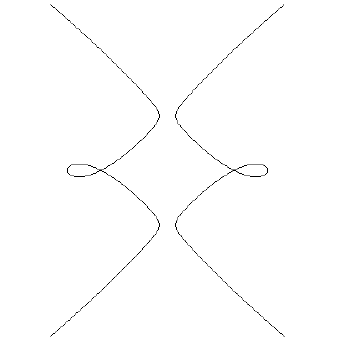}
  \end{center}
\end{figure}

The reciprocal polar, with respect to the standard flag and standard quadric,
$S_{3}^{\prime
\prime}=\mathcal{V}(X_2\frac{
\partial f_{3}}{
\partial X_1}-X_1\frac{
\partial f_{3}}{
\partial X_2})$, intersects each component in non-singular points, since the
double points are not among the points on each component with the locally
shortest or longest distance from the origin.

\begin{figure}[htbp!]
\caption{The curve $S_{3}$ and its reciprocal polar $S_{3}^{\prime \prime}$.}
  \label{ikkekomp2dp_polar}
  \begin{center}
  \includegraphics{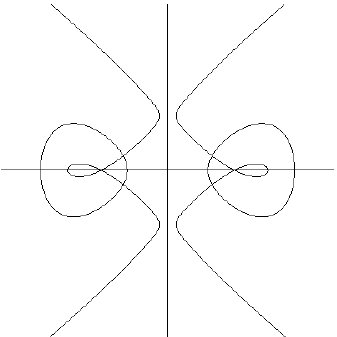}
  \end{center}
\end{figure}
\medskip

{\sc Example 4.} This example illustrates Proposition \ref{transform}. We
start with the irreducible affine curve $S_{4}$ defined by the polynomial
\begin{small}
\[
f_{4}:=X_1^2-X_2(X_2+1)(X_2+2).
\]
\end{small}
 This curve passes through the origin, so the
reciprocal polar variety, when $Q$ is the standard quadric and
$L=\mathcal{V}(X_0)$ is the line at infinity, will not contain any points
from the connected component containing the origin, since we are not counting
points which are on both the variety and the flag. Instead we choose the point
$(1,0)$ (or $(1:1:0)$ in projective coordinates), and we must find a polynomial
$q^{\prime}$ such that $(1:1:0)^{\perp ^{\prime}}$ is the line $\mathcal{V}(X_0)$. We see
that the polynomial $2X_0^2-2X_0X_1+X_1^2+X_2^2$ will do, and we will now
find the reciprocal polar variety $W_{L}^{\perp ^{\prime}}(S_{4})$. If 
$P=(1:p_1:p_2)$, the point $\langle P,
L^{\perp ^{\prime}} \rangle^{\perp ^{\prime}}=P^{\perp ^{\prime}}\cap L$ is the point $(0:p_2:1-p_1)$,
so the affine cone over it, $I_{P,L^{\perp ^{\prime}}}$, is the
line $\mathcal{V}((p_1-1)X_1+p_2X_2)$. 
The reciprocal polar variety $W_{L}^{\perp ^{\prime}}(S_{4})$ is the set
\begin{displaymath}
\textstyle
\{ P \in (S_{4})_{\mathbb{R}} \,|\, p_2\frac{
\partial f_{4}}{
\partial X_1}(P)+(1-p_1)\frac{
\partial f_{4}}{
\partial X_2}(P)=0\};
\end{displaymath} 
this set consists of the points on $(S_{4})_{\mathbb{R}}$ with
the locally shortest or longest Euclidean distance to the point $(1,0)$, and it
contains points from each component of $(S_{4})_{\mathbb{R}}$.

\begin{figure}[htbp!]
\caption{The curve $S_{4}$.}
  \label{transformkurve}
  \begin{center}
    \includegraphics{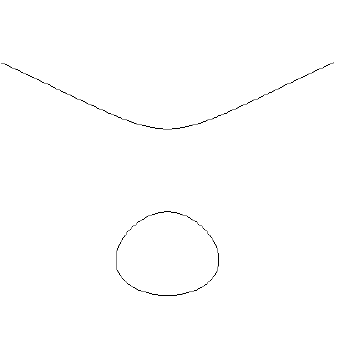}
  \end{center}
\end{figure}

\begin{figure}[htbp!]
\caption{The curve $S_{4}$ and its reciprocal polar $S_{4}^{\prime 
\prime}$.}
 \label{transformpolar}
  \begin{center}
    \includegraphics{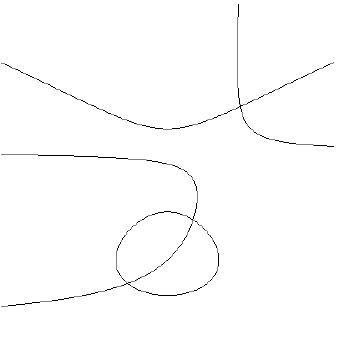}
  \end{center}
\end{figure}
\medskip

{\sc Example 5.} This is an example showing that Proposition
\ref{proprpv} does not hold for curves with arbitrary singularities.
Consider the curve  $S_{5}$ given by the
polynomial 
\begin{small}
\begin{displaymath} 
\textstyle
f_{5}=((X_1-4)^2+(X_2-2)^2-1)^2+\frac{1}{100}((X_1-\frac{7}2)(X_1-\frac{9}2))^3, 
\end{displaymath}
\end{small}
The real components of this curve do not contain points on
the reciprocal polar variety other than the four cusps. This can be
seen by calculating the intersection points between $S_{5}$ and its reciprocal
polar
$S_{5}^{\prime \prime}=\mathcal{V}(X_1\frac{\partial{f_{5}}}{\partial X_2}
-X_2\frac{\partial{f_{5}}}{\partial
X_1})$.
 
 \begin{figure}[htbp!]
  \caption{The curve $S_{5}$.}
  \label{dualpolareks}
  \begin{center}
        \includegraphics{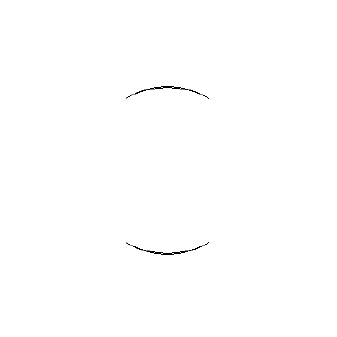}
  \end{center}
\end{figure}

\begin{figure}[htbp!]
  \caption{The reciprocal polar $S_{5}^{\prime \prime}$.}
  \label{kurveogdualpolar}
  \begin{center}
    \includegraphics{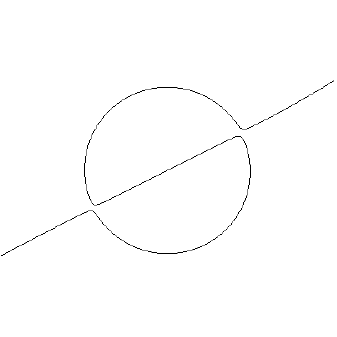}
  \end{center}
\end{figure}
\bigskip

\bibliography{PolarRef}
\bibliographystyle{siam}

\end{document}